
\def\Z_2{{\bf Z}_2}
\def\Z{{\bf Z}}
\def\C{{\bf C}}
\def\>{\rangle}

\baselineskip=16pt
\leftskip 36pt

\noindent
{\bf UNITARIZABLE REPRESENTATIONS OF THE DEFORMED \hfill\break
PARA-BOSE SUPERALGEBRA  U$_q$[osp(1/2)] AT ROOTS OF 1}

\vskip 32pt
\noindent
Short title: {\bf UNITARIZABLE ROOT OF 1 IRREPS OF U$_q$[osp(1/2)]}

\vskip 32pt
\noindent
T. D. Palev{*} and N. I. Stoilova{*}

\noindent
International Centre for Theoretical Physics, 34100 Trieste, Italy
\vskip 12pt

\footnote{*}{Permanent address: Institute for Nuclear Research and
Nuclear Energy, 1784 Sofia, Bulgaria; E-mail:
palev@bgearn.bitnet, stoilova@bgearn.bitnet}

\vskip 4in
\noindent
Classification numbers according to the  Physics and Astronomy
Classification Scheme: 02.10.Tr, 02.20.Fh, 03.65.Fd.

\vfill \eject

\vskip 48pt

{\bf Abstract.} The unitarizable irreps of the deformed para-Bose
superalgebra $pB_q$, which is isomorphic to $U_q[osp(1/2)]$,
are classified at $q$ being root of 1. New
finite-dimensional irreps of $U_q[osp(1/2)]$ are found. Explicit
expressions for the matrix elements are written down.

\vfill \eject

\vskip 48pt

\noindent
{\bf 1. Introduction}
\bigskip

\noindent
In the present paper we study unitarizable root of unity
representations of the Hopf algebra $pB_q(1)\equiv pB_q$,
introduced in [1]. It is generated (essentially) by one pair of
deformed para-Bose  operators $a^\pm$. The irreps of
$pB_q$ at generic values of $q$ are infinite-dimensional and are
realized in deformed para-Bose Fock spaces $F(p), \; p\in {\bf
C}$ [1]. The multimode Hopf algebra $pB_q(n)$, corresponding to
$n$ pairs of deformed pB operators $a_1^\pm,
a_2^\pm,\ldots,a_n^\pm$ was defined in [2-5]. The case of any
number of deformed para-Fermi operators was worked out in [6].

So far various deformations of  para-Bose and para-Fermi
statistics were considered from different points of view [7-22].
Some of them are not related to any Hopf algebra structure.
A guiding principle of the approaches in [1-5], which we follow,
is to preserve, similar to the nondeformed case [23], the
identification of $pB_q(n)$ with $U_q[osp(1/2n)]$: $pB_q(n)$ is
an associative superalgebra isomorphic (as a Hopf algebra) to the
deformed universal enveloping algebra  $U_q[osp(1/2n)]$ of the
orthosymplectic Lie superalgebra $osp(1/2n)$.

The Hopf algebra structure of $pB_q(n)$ has an important
advantage: using the comultiplication, one can define new
representations of the deformed operators (and hence of
$U_q[osp(1/2n)]$) in any tenzor product of representation spaces.
In particular one can use the Fock space of $n$ pairs of
commuting deformed Bose operators [24-27], since they give a
representation of $U_q[osp(1/2n)]$ [28]. Even in the nondeformed
case the only effective technique for constructing
representations of parabosons or of $osp(1/2n)$ (for large $n$)
is through tenzor products of bosonic Fock spaces (see [29] for
more disscusions in this respect).

The definition of $U_q[osp(1/2n)]$ in terms of
its Chevalley generators is well known [30-35].
Although for $n>1$ the deformed pB operators
$a_1^\pm, a_2^\pm,\ldots,a_n^\pm$
are very different from the Chevalley generators, the relations
determining $U_q[osp(1/2n)]$ through $a_1^\pm, a_2^\pm,\ldots,a_n^\pm$
are not more involved [4,5].
At $n=1$, namely in the case we consider, $a^\pm$ are
proportional to the Chevalley generators of $pB_q=U_q[osp(1/2)]$.

The finite-dimensional irreps of $U_q[osp(1/2)]$ at generic $q$
were constructed in [36, 37].  Some root of unity highest weight
irreps were also obtained in [37]; both highest weight and cyclic
representations were studied in detail in [38-40].  

Our carrier representation spaces $F(p), \; p\in {\bf C}$ will be
deformed Fock spaces [1], which are in fact the Verma
modules used in [39]. At root of unity cases each such space is no
more irreducible; it contains infinitely many invariant
subspaces. The irreps are realized  in appropriate factor spaces
of $F(p)$ with the vacuum being the highest weight vector.

To our best knowledge the root of 1 irreps of $U_q[osp(1/2)]$
obtained in Section 3.2 and those labeled with an integer $p$
[Sections  3.1{\it a}, 3.1{\it b}, 3.3] were not described in the
literature so far.  Our other main result is the classification
of  the unitarizable Fock irreps of $pB_q$ (= unitarizable  Verma
representations of $U_q[osp(1/2)]$) at roots of 1 (see (4.7)). We
write down explicit expressions for the transformation of the
basis under the action of the deformed pB generators.

The reason to pay a special attention to the unitarizable
representations stems from physical considerations. In all known
to us applications of deformed parastatistics [7-22] it is
assumed that the Hermitian conjugate $(a^-)^\dagger$ of the annihilation
operator $a^-$ equals to the creation operator $a^+$,
$$
(a^-)^\dagger=a^+. \eqno(1.1)
$$
In case of deformed para-oscillators [9, 15], for instance, or
more generally in any deformed quantum mechanics (see for
instance [41] and the references therein) the unitarity condition
(1.1) is equivalent (as in the canonical case) to the
requirement the position and the momentum operators to be
selfadjoint operators.

The paper is organized as follows. In Sec. 2 we recall the
definition of the deformed para-Bose algebra and its Fock
representations at generic $q$. Sec. 3 is devoted to a detailed
study of the root of 1 irreps. The indecomposible representations
both finite-dimensional and infinite-dimensional are also
mentioned. The unitarizable representations are classified in
Sec. 4. Sec. 5 contains some concluding remarks.

Throughout we use the following abbreviations and notation:

\vskip 6pt

${\bf C}$ --- all complex numbers

${\bf Z}$ --- all integers

${\bf Z}_+$ --- all nonnegative integers

${\bf Z}_2 = \{\bar{0},\bar{1}\}$ --- the ring of all integers modulo 2

$[A,B]=AB-BA, \hskip 6pt \{A,B\}=AB+BA$.

\vskip 32pt
\noindent
{\bf 2. The para-Bose Hopf algebra pB$_q$ and its Fock representations}

\bigskip
\noindent
To begin with we  summarize some of the results from [1],
changing slightly the notation.

\smallskip
\noindent
{\it Definition 1.} The para-Bose algebra $pB_q,\;
q\in{\bf C} \backslash \{0, \pm 1\}$, is the associative superalgebra
over {\bf C} with unit 1 defined by the following generators and
relarions
$$
\eqalignno{
& {\rm Generators}: a^\pm,\; K^{\pm 1} & (2.1)  \cr
& {\rm Relations}: KK^{-1}=K^{-1}K=1\quad
  Ka^\pm=q^{\pm 2}a^\pm K  \quad
  \{a^+,a^-\}={{K-K^{-1}}\over{q-q^{-1}}} & (2.2)  \cr
& {\bf Z}_2-{\rm grading}: deg(K^{\pm 1})=\bar{0} \quad
  deg(a^\pm)=\bar{1}. & (2.3)    \cr
}
$$
\bigskip
The creation operator $a^+$ (resp. the annihilation operator
$a^-$) is the negative (resp. the positive) root vector of
$pB_q=U_q[osp(1/2)]$.

Setting $K=q^H$ with $q=e^\eta, \; \eta \in {\bf C}$, one recovers
as $q\rightarrow 1 \;\;(\eta \rightarrow 0)$ the defining relations
of the nondeformed para-Bose operators [42]
($\xi, \eta, \epsilon =\pm \; {\rm or} \; \pm 1 $):
$$
[\{a^\xi,a^\eta\},a^\epsilon]=(\epsilon - \eta)a^\xi
+(\epsilon - \xi)a^\eta \eqno(2.4)
$$

\noindent
with $H=\{a^+,a^-\}$.

It was already shown in [1] how $pB_q$ can
be endowed with a comultiplication, a counity and an antipode;
here we shall be not concerned with this additional structure.

The (deformed) Fock space $F(p)$ is defined for any complex
number $p$, $p \in {\bf C}$,
postulating that
$F(p)$ contains a vacuum vector $\vert p;0\rangle$, namely
$$
a^-\vert p;0\rangle=0
\eqno(2.5)
$$

\noindent
so that
$$
K\vert p;0\rangle=q^p\vert p;0\rangle \quad
(\Leftrightarrow  H\vert p;0\rangle=p\vert p;0\rangle). \eqno(2.6)
$$

\noindent
At the limit $q \rightarrow 1$ the above definition of $F(p)$
reduces to the usual one   \hfill\break
($a^-\vert p;0\rangle=0, \; a^-a^+\vert
p;0\rangle=p\vert p;0\rangle)$, where $p$ is the order of the
parastatistics [42].

$F(p)$ is an infinite-dimensional linear space with a basis
$$
\vert p;n\rangle=(a^+)^n\vert p;0\rangle
\quad n\in {\bf Z}_+
\eqno(2.7)
$$
\noindent
and a highest weight vector $\vert p;0\rangle$.

Setting
$$
[n]={{q^n-q^{-n}}\over{q-q^{-1}}} \quad
\{n\}={{q^n+q^{-n}}\over{q+q^{-1}}} \eqno(2.8)
$$
\noindent
one can write the transformation of the basis as follows:
$$
\eqalignno{
& K\vert p;n\rangle=q^{2n+p}\vert p;n\rangle & (2.9a) \cr
& a^+\vert p;n\rangle=\vert p;n+1\rangle & (2.9b) \cr
& a^-\vert p;n\rangle=[n]\{n+p-1\}\vert p;n-1\rangle
  \quad {\rm for}\;n={\rm even \; number} & (2.9c)\cr
& a^-\vert p;n\rangle=[n+p-1]\{n\}\vert p;n-1\rangle
  \quad {\rm for}\;n={\rm odd \; number}. & (2.9d)\cr
}
$$
At generic $q$ each Fock space $F(p),\; p\in {\bf C}$, is an
infinite-dimensional simple (=irreducible) $pB_q$ module [1].

\vskip 32pt
\noindent
{\bf 3. Root of unity representations}

\bigskip
\noindent
If $q$ is a root of $1$ the $pB_q$ module $F(p)$ may no longer be
irreducible. More precisely,
\smallskip
\noindent
{\it Proposition 1.} The Fock space $F(p)$ is indecomposible if
and only if $q=e^{i{\pi \over 2}{m\over k}}$ for every
$m,k \in {\bf Z}$ such that  $q \notin \{\pm1 , \pm i\}$, i.e., $m\neq
0(mod\; k)$.

\smallskip
\noindent
{\it Proof.} We exclude from consideration $q \in \{ \pm 1, \pm i
\}$ since at those values of $q$ the expresions $(2.9)$ are not
defined (at $q=\pm 1$ also $pB_q$ is undefined).

If $q\neq e^{i{\pi \over 2}{m\over k}}$ the coefficients in front
of $\vert p;n\rangle$ in (2.9{\it c,d}) vanish only for $n=0$. Hence the only
singular vector is the vacuum, i.e., $F(p)$ is an irreducible module.
If $q=e^{i{\pi \over 2}{m\over k}}$ then, for instance, $\vert
p;2k\rangle$ is a singular vector, $a^-\vert p;2k\rangle=0$. Therefore the
proper subspace of $F(p)$, spanned on $\vert p;n\rangle, \; n\geq 2k$ is
an invariant subspace. Then eq.  (2.9{\it b}) yields that the
representation is indecomposible. This completes the proof.

\smallskip
The algebras $pB_q$ corresponding to all possible values of $m$
and $k$ contain several isomorphic copies. Clearly we can always
assume that $k>0$ and that $m$ and $k$ are co-prime, i.e.,
$m\over k$ is an irreducible fraction. Further we note that
the algebras $pB_q$ and $pB_{\tilde q^{\xi}}$ are isomorphic for
$$
q=e^{i{\pi \over 2}{m\over k}}\; {\rm and}\;
\tilde q^\xi =e^{i{\pi \over 2}{{2k+\xi m}\over k}} \eqno(3.1)
$$
whenever $\xi =+$ and $m=1,2, \ldots , 2k-1$ or $\xi =-$ and
$m=1, 2,\ldots, k-1$, since the generators ${\tilde a}^\pm=a^{\pm
\xi}$ and $\tilde K = -\xi K$ of $pB_{\tilde q^{\xi}}$  satisfy
the defining relations (2.1)-(2.3) for $pB_q$. The case $\xi=+$
indicates that we can set $m \in
\{1,2,\ldots,2k-1\}$; the case $\xi=-$ further shows that
without loss of generality we can assume that
$m \in \{1,2,\ldots,k-1\}$. The case $k=1$ is excluded from these
conditions. Thus, without losing any of the algebras $pB_q$ for
which the Fock space $F(p)$ is indecomposible, we restrict
$m$ and $k$ to values, which we call admissible. The
fraction $m\over k$ is said to be admissible if
$$
\eqalignno{
& k=2,3,\ldots  & (3.2a)\cr
& {m\over k} \in  \Big\{ {1\over k},
  {2\over k}, \ldots, {{k-1}\over k} \Big\}  & (3.2b) \cr
& {m\over k} \; {\rm is \; an \; irreducible  \; fraction} \;
  (m \; {\rm and} \;  k \; {\rm are \; co\!-\!prime)}. & (3.2c) \cr
}
$$

Passing to a discussion of the root of 1
representations, we first note that the
vectors
$$
\vert p;0\rangle,\vert p;2k\rangle, \vert p;4k\rangle,\ldots,
\vert p;2kN\rangle, \ldots   \eqno(3.3)
$$
are singular vectors in $F(p)$. Indeed $[2kN]=0$ and therefore
(see (2.9c))
$$
a^- \vert p;2kN\rangle=0 \quad N\in {\bf Z}_+. \eqno(3.4)
$$
\noindent
The subspaces
$$
V_{\vert p;2kN\rangle}=span\{\vert p;n\rangle \vert n\geq 2kN\}
\eqno(3.5)
$$
\noindent
are infinite-dimensional invariant subspaces of $F(p)$ with
highest weight vectors $\vert p; 2kN\rangle$. Clearly
$$
F(p)=V_{\vert p;0\rangle}\supset V_{\vert p;2k\rangle}\supset
V_{\vert p;4k\rangle}\supset
\ldots \supset V_{\vert p;2kN\rangle} \supset \ldots .  \eqno(3.6)
$$
For each $N<M \in {\bf Z}_+$
define a $2(M-N)k-$dimensional factor space
$$
W_{\vert p;2kN\rangle,M}=V_{\vert p;2kN\rangle}/V_{\vert p;2kM\rangle}.
\eqno(3.7)
$$
\noindent
Let $\xi_x$ be the equivalence class of $x\in \xi_x$. The vectors
$$
\xi_{\vert p;2kN\rangle},\; \xi_{\vert p;2kN+1\rangle},\;
\xi_{\vert p;2kN+2\rangle},
\ldots ,\; \xi_{\vert p;2kM-1\rangle} \eqno(3.8)
$$
constitute a basis in $W_{\vert p;2kN\rangle,M}$. The relations
$$
a^\pm \xi_{\vert p;n\rangle}=\xi_{a^\pm \vert p;n\rangle} \quad
K \xi_{\vert p;n\rangle}=\xi_{K \vert p;n\rangle} \eqno(3.9)
$$
endow $W_{\vert p;2kN\rangle,M}$ with a structure of a $pB_q$
module. Observe that $a^+ \xi_{\vert p;2kM-1\rangle}=0$.

We  shall simplify the notation, identifying $\xi_{\vert p;n\rangle}$
with its representative $\vert p;n\rangle$. Then, adding to eqs. (2.9)
the condition $a^+\vert p;2kM-1\rangle=0$, one obtains the
transformations of $W_{\vert p;2kN\rangle,M}$.

We summarise.  The space $W_{\vert p;2kN\rangle,M}, \; N<M \in {\bf
Z}_+$ is a $2(M-N)k$-dimensional $pB_q$ module with a highest
weight $\vert p;2kN\rangle$, a basis
$$
\vert p;n> \quad n=2kN, 2kN+1,\ldots,2kM-1 \eqno(3.10)
$$
and transformation relations
$$
\eqalignno{
& K\vert p;n\rangle=q^{2n+p}\vert p;n\rangle & (3.11a) \cr
& a^+\vert p;n\rangle=\vert p;n+1\rangle \;\; n\neq 2kM-1 & (3.11b) \cr
& a^+\vert p;2kM-1\rangle=0 & (3.11c) \cr
& a^-\vert p;n\rangle=[n]\{n+p-1\}\vert p;n-1\rangle
  \quad {\rm for}\;n={\rm even \; number} & (3.11d)\cr
& a^-\vert p;n\rangle=[n+p-1]\{n\}\vert p;n-1\rangle
  \quad {\rm for}\;n={\rm odd \; number}. & (3.11e)\cr
}
$$

\noindent
If $M-N>1$, eqs. (3.11) define an indecomposible
finite-dimensional representation of $pB_q$ in
$W_{\vert p;2kN\rangle,M}$ for any $p\in {\bf C}$. If
$M-N=1$, $W_{\vert p;2kN\rangle,N+1}$ is either irreducible or
indecomposible. In this case the vectors
$$
\vert p;2kN\rangle,\;\vert p;2kN+1\rangle,\;\vert p;2kN+2\rangle,\ldots,
\vert p;2k(N+1)-1\rangle \eqno(3.12)
$$
constitute a basis in $W_{\vert p;2kN\rangle,N+1}$.

\smallskip
\noindent
{\it Proposition 2.} For each $s\in {\bf Z}$ and $N\in {\bf Z}_+$
the $pB_q$ modules $W_{\vert p;0\rangle,1}$ and
$W_{\vert p+4ks;2kN\rangle,N+1}$ are equivalent;
they carry one and the same  $2k-$dimensional representation of
$pB_q$.

\noindent
{\it Proof.} Define a one-to-one map
$$
\varphi(\vert p;n\rangle)=\vert p+4ks;n+2kN\rangle
\quad n=0,1,2,\ldots,2k-1 \eqno(3.13)
$$
of the basis in $W_{\vert p;0\rangle,1}$ onto the basis of
$W_{\vert p+4ks;2kN\rangle,N+1}$
and extend it to a linear map from $W_{\vert p;0\rangle,1}$
onto $W_{\vert p+4ks;2kN\rangle,N+1}$.
{}From (3.11) one derives that $\varphi$ commutes with $pB_q$,
$$
a\varphi(\vert p;n\rangle)=\varphi(a\vert p;n\rangle)
\quad a=a^\pm,\; K. \eqno(3.14)
$$
i.e., $\varphi$ is an intertwining operator.
In particular the matrices of the generators $a^\pm$ and $K$ are
the same in the basis
$$\vert p;0\rangle,\vert p;1\rangle,
\vert p;2\rangle,\ldots,\vert p;2k-1\rangle \eqno(3.15)
$$
of $W_{\vert p;0\rangle,1}$ and in the basis
$\vert p+4ks;2kN\rangle, \vert p+4ks;2kN+1\rangle, \ldots,\vert
p+4ks;2k(N+1)-1\rangle$ of $W_{\vert p+4ks;2kN\rangle,N+1}$,
respectively. This completes the proof.

In view of Proposition 2 we shall consider from now on only the
vacuum modules $W_{\vert p;0\rangle,1}$, restricting also the
values of $p$ to the interval
$$
0<Re(p)\leq 4k. \eqno(3.16)
$$
\noindent
We  often write $W(\vert p;m\rangle,\vert p;n\rangle)$ whenever we
wish to indicate that
$$
\vert p;m\rangle,\vert p;m+1\rangle,\vert p;m+2\rangle,
\ldots,\vert p;n\rangle \eqno(3.17)
$$
is a basis in the linear space $W(\vert p;m\rangle,\vert p;n\rangle)$.
In view of (3.15)
$$
W_{\vert p;0\rangle,1}=W(\vert p;0\rangle,\vert p;2k-1\rangle). \eqno(3.18)
$$

We proceed to study in detail the structure of the Fock spaces
for different admissible values of  $m\over k$.
To this end we will consider three cases:
3.1  $k$ is even, $m$ is odd; 3.2  $k$ is odd, $m$ is odd;
3.3 $k$ is odd, $m$ is even.

\vskip 16pt
\noindent
{\bf 3.1 The case k=even, m=odd}
\bigskip
\noindent
If $p$ is not an integer the coefficients $\{n+p-1\}$ in $(2.9c)$
and $[n+p-1]$ in $(2.9d)$ never vanish. Therefore the vectors
(3.3) are the only singular vectors in $F(p)$ and the vacuum
$\vert p;0\rangle$ is the only singular vector in $W_{\vert p;0\rangle,1}$.
The eigenvalues of $K$ on $\vert p;0\rangle$ are different for
different $p$, obeying (3.16). This gives rise to the following

\smallskip
\noindent
{\it Proposition 3.} The $pB_q$ modules $W(\vert p;0\rangle,\vert
p;2k-1\rangle)$ are simple  for $p\notin \{1,2,\ldots,4k\}$. All of
them are $2k-$dimensional. The irreps corresponding to different
$p$ from (3.16) are inequivalent.

The transformation of the basis (3.15)
follows from eqs. (3.11) at $N=M-1=0$. The relations bellow are
written in a slightly more general form in order to accomodate
also other cases. For
$$
\vert p;0\rangle,\;\vert p;1\rangle,\ldots, \vert p;L\rangle
\eqno(3.19a)
$$
set
$$
\eqalignno{
& K\vert p;n\rangle=q^{2n+p}\vert p;n\rangle & (3.19b) \cr
& a^-\vert p;n\rangle=[n]\{n+p-1\}\vert p;n-1\rangle
  \quad {\rm for}\;n={\rm even \; number} & (3.19c)\cr
& a^-\vert p;n\rangle=[n+p-1]\{n\}\vert p;n-1\rangle
  \quad {\rm for}\;n={\rm odd \; number}. & (3.19d)\cr
& a^+\vert p;n\rangle=\vert p;n+1\rangle \;\; n\neq L & (3.19e) \cr
& a^+\vert p;L\rangle=0. & (3.19f) \cr
}
$$
If $L=2k-1$ the above equations give the transformation of
$W(\vert p;0\rangle,\vert p;2k-1\rangle)$.

\vskip 14pt
\noindent
{\it 3.1a Representations with even $p$ }
\bigskip
\noindent
All modules $W(\vert p;0\rangle \vert p;2k-1\rangle)$ corresponding to even
values of $p$ are no more irreducible.
\smallskip
\noindent
{\it Proposition 4.} To each $p\in \{2,4,\ldots,4k\}$ there
corresponds a simple $pB_q $ module
$W(\vert p;0\rangle),\vert p;L\rangle)$ with a basis $\vert
p;0\rangle,\vert p;1\rangle,\ldots,\vert p;L\rangle$ and values
of $L$ as follows:
$$
\eqalignno{
& L=2k-p \quad {\rm for}
  \quad p\in \{2,4,\ldots,2k\} & (3.20)\cr
& L=4k-p \quad {\rm for}
\quad p\in \{2k+2,2k+4,\ldots,4k\}. & (3.21)\cr
}
$$
\noindent
The transformation of the basis is described with the  eqs.
(3.19) for the above values of $L$. All $2k$ modules
carry different, inequivalent irreps of $pB_q$.
\smallskip
\noindent
{\it Proof.} We consider in detail the case
$p\in \{2,4,\ldots,2k\}$. The  module $W(\vert p;0\rangle,\vert
p;2k-1\rangle)$ contains only two singular
vectors $\vert p;0\rangle$ and $\vert p;2k-p+1\rangle$:
$$
a^-\vert p;0\rangle=0 \quad a^-\vert p;2k-p+1\rangle=0. \eqno(3.22)
$$
In view of $(3.19e)$ $W(\vert p;0\rangle,\vert p;2k-1\rangle)$ is
indecomposible. Its invariant subspace  \hfill\break
$W(\vert p;2k-p+1\rangle,\vert p;2k-1\rangle)$
is simple. The factor space
$$
W(\vert p;0\rangle,\vert p;2k-1\rangle)/
W(\vert p;2k-p+1\rangle,\vert p;2k-1\rangle)
$$
with a basis
$$
\xi_{\vert p;0\rangle},\xi_{\vert p;1\rangle},
\ldots,\xi_{\vert p;2k-p\rangle} \eqno(3.23)
$$
is turned into an irreducible $pB_q$ module, setting
$$
a\xi_{\vert p;n\rangle}=\xi_{a\vert p;n\rangle} \quad a=a^\pm,\;K. \eqno(3.24)
$$
Therefore
$$
a^+{\xi_{\vert p;2k-p\rangle}}=\xi_{\vert p;2k-p+1\rangle}=0. \eqno(3.25)
$$
\noindent
As before we identify the equivalence classes with their
representatives and in particular
$$
\xi_{\vert p;n\rangle}=\vert p;n\rangle. \eqno(3.26)
$$
Then
$$
W(\vert p;0\rangle,\vert p;2k-1\rangle)/
W(\vert p;2k-p+1\rangle,\vert p;2k-1\rangle)=
W(\vert p;0\rangle,\vert p;2k-p\rangle). \eqno(3.27)
$$

{}From (3.24) and (3.25) one derives the transformation relations
of $W(\vert p;0\rangle,\vert p;2k-p\rangle)$, which are given with eqs.
(3.19) for $L=2k-p$. The transformations of the invariant subspaces
are described also with (3.19), but for $n=2k-p+1, 2k-p+2,\ldots,
2k-1=L$.

The cases with $p\in \{2k+2, 2k+4,\ldots, 4k\}$ are similar. The
only singular vectors in \hfill\break
$W(\vert p;0\rangle,\vert p;2k-1\rangle)$ are the vacuum
$\vert p;0\rangle$ and $\vert p;4k-p+1\rangle$. Therefore the invariant
subspace $W(\vert p;4k-p+1\rangle,\vert p;2k-1\rangle)$ is
irreducible and its transformations are described with (3.19) for
$n=4k-p+1, 4k-p+2,\ldots, 2k-1=L$. The factor space $W(\vert
p;0\rangle,\vert p;4k-p\rangle)$ is also irreducible and transforms
according to (3.19) with $L=4k-p$.

So far we have four kinds of simple $pB_q$ modules.
We proceed to show that
$$
\eqalignno{
& W(\vert p;2k-p+1\rangle,\vert p;2k-1\rangle)=
  W(\vert p';0\rangle,\vert p';4k-p'\rangle) & \cr
& \hskip 12mm {\rm for} \quad  p=4k-p'+2\in \{2,4,\ldots,2k\} & (3.28)\cr
& W(\vert p;4k-p+1\rangle,\vert p;2k-1\rangle)=
  W(\vert p';0\rangle,\vert p';2k-p'\rangle) & \cr
& \hskip 12mm {\rm for} \quad  p=4k-p'+2\in
  \{2k+2,2k+4,\ldots,4k\} & (3.29)\cr
}
$$
\noindent
where the equality means that the corresponding modules are
equivalent, they carry equivalent representations of $pB_q$.

In order to prove (3.28) we set a one to one linear map
$$
\varphi:\;\; W(\vert p;2k-p+1\rangle,\vert p;2k-1\rangle)\; \rightarrow \;
  W(\vert p';0\rangle,\vert p';4k-p'\rangle) \eqno(3.30)
$$
\noindent
which restricted to the corresponding basisses read:
$$
\varphi(\vert p;n\rangle)=\vert p';n'\rangle, \quad n=2k-p+1,2k-p+2,
\ldots 2k-1 \Longleftrightarrow n'=0,1,2,\ldots,4k-p'. \eqno(3.31)
$$
\noindent
It is straighforward to check, using eqs. (3.19), that $\varphi$
commutes with the $pB_q$ generators, i.e., it is an intertwining operator.

The proof of (3.29) is similar. It remains to show that all $2k$
modules (3.20) and (3.21), i.e.,
$$
W(\vert p;0\rangle,\vert p;2k-p\rangle)\;\; {\rm and}\; \;
W(\vert p+2k;0\rangle,\vert p+2k;2k-p\rangle)
\quad p\in \{2,4,\ldots ,2k\} \eqno(3.32)
$$
are inequivalent. To this end we observe that the  modules
corresponding to different $p$ in (3.32) have different
dimensions. The modules with the same dimensions, namely
$W(\vert p;0\rangle,\vert p;2k-p\rangle)$ and $W(\vert p+2k;0\rangle,\vert
p+2k;2k-p\rangle)$ have different spectrum of the Cartan generator $K$
and therefore are also inequivalent. This completes the proof.


It has been noted in [37, 43] that the Casimir operator of
$U_q[osp(1/2)]$ is no longer sufficient to label the root of
unity representations. In particular this is the case with the
modules (3.32), which have the same dimension. The Casimir
operator [36] reads in our notation:
$$
\eqalignno{
2C_2(q) & = q^2K^2+q^{-2}K^{-2} & \cr
& +(q^2-q^{-2})(q-q^{-1})(q^2K+q^{-2}K^{-1})a^-a^+
-(q^2-q^{-2})^2(a^-)^2(a^+)^2. &(3.33)\cr
}
$$
Its eigenvalue is one and the same,
$$
C_2(q)={1\over 2}(q^{2p-2}+q^{-2p+2})+2 \eqno(3.34)
$$
on four inequivalent modules, namely

$$
\eqalign{
& \hskip 6pt W(\vert p;0\rangle,\vert p;2k-p\rangle)\;\; \;\;
   W(\vert p+2k;0\rangle,\vert p+2k;2k-p\rangle) \quad dim=2k-p+1 \cr
& \hskip 6pt W(\vert 2k-p+2;0\rangle,\vert 2k-p+2;p-2\rangle)\;\; \;\;
   W(\vert 4k-p+2;0\rangle,\vert 4k-p+2;p-2\rangle) \quad dim=p-1. \cr
}\eqno(3.35)
$$

\noindent
Observe that pairwise the dimensions of these modules are different.

Let us add that the additional central elements [37, 38] $(a^\pm)^{4k}$
and $(K)^{2k}$ do not distinguish among the inequivalent modules
with the same dimension. In fact the operators $(a^\pm)^{4k}$ vanish
within each simple module $W(\vert p;0 \rangle,\vert p;L \rangle)$.

\vfill\eject
\noindent
{\it 3.1b Representations with odd $p$ }
\bigskip
\noindent
{\it Proposition 5.} To each $p\in \{1,3,\ldots,4k-1\}$ there
corresponds an irreducible $pB_q $ module \hfill\break
$W(\vert p;0\rangle,\vert p;L\rangle)$ with a basis $\vert
p;0\rangle,\vert p;1\rangle,\ldots,\vert p;L\rangle$ and values
of $L$ as follows:

$$
\eqalignno
{
& L=k-p \quad \;\; {\rm for}
\quad  p\in\{1,3,\ldots,k-1\} & (3.36a)\cr
& L=3k-p \quad {\rm for}
\quad p\in\{k+3,k+5,\ldots,3k-1\} & (3.36b)\cr
& L=5k-p \quad {\rm for}
\quad p\in\{3k+3,3k+5,\ldots,4k-1\} & (3.36c) \cr
& L=2k-1 \quad {\rm for}
\quad p=k+1,\;3k+1. &  (3.36d)\cr
}
$$
\noindent
The transformations of the basis are given with the  eqs.
(3.19) for the corresponding values of $L$, indicated above. All
these $2k$ modules carry inequivalent irreps of $pB_q$.
\smallskip
\smallskip
\noindent
{\it Proof.} The proof is similar to the one of Proposition 4.
We stress on certain points only. In this case we have both simple
and indecomposible modules $W(\vert p;0\rangle,\vert p;2k-1\rangle)$.
More precisely,
$$
W(\vert p;0\rangle,\vert p;2k-1\rangle)  \quad \quad
\quad p=k+1,\;3k+1   \eqno(3.37)
$$
are irreducible; the rest of the modules
$$
W(\vert p;0\rangle,\vert p;2k-1\rangle) \quad p\in\{1,3,5,\ldots,4k-1\}
\;\;p\neq k+1,\;3k+1 \eqno(3.38)
$$
are indecomposible. Each module in (3.38) contains apart from the
vacuum only one more singular vector. Therefore the invariant
subspaces and the factor spaces are irreducible. Also here each
invariant subspace is equivalent to a factor space. More precisely,
$$
\eqalignno{
& W(\vert p;k-p+1\rangle,\vert p;2k-1\rangle)=
W(\vert p';0\rangle,\vert p';3k-p'\rangle)&\cr
& \hskip 12mm {\rm for} \quad  p=2k-p'+2\in
\{1,3,\ldots,k-1\} & (3.39) \cr
}
$$
$$
\eqalignno{
& W(\vert p;3k-p+1\rangle,\vert p;2k-1\rangle)=
W(\vert p';0\rangle,\vert p';k-p'\rangle)&\cr
& \hskip 12mm {\rm for} \quad  p=2k-p'+2\in
\{k+3,k+5,\ldots,2k+1\} & (3.40) \cr
}
$$
$$
\eqalignno{
& W(\vert p;3k-p+1\rangle,\vert p;2k-1\rangle)=
W(\vert p';0\rangle,\vert p';5k-p'\rangle)&\cr
& \hskip 12mm {\rm for} \quad  p=6k-p'+2\in
\{2k+3,2k+5,\ldots,3k-1\} & (3.41) \cr
}
$$
$$
\eqalignno{
& W(\vert p;5k-p+1\rangle,\vert p;2k-1\rangle)=
W(\vert p';0\rangle,\vert p';3k-p'\rangle)&\cr
& \hskip 12mm {\rm for} \quad  p=6k-p'+2\in
\{3k+3,3k+5,\ldots,4k-1\}. & (3.42) \cr
}
$$
Therefore (up to equivalence) we are left only with the vacuum
modules (3.36). Among them one finds all of the time pairs with
the same dimension, namely
$$
\eqalignno{
& dim[W(\vert p;0\rangle,\vert p;k-p\rangle)]=
dim[W(\vert p+2k;0\rangle,\vert p+2k;k-p\rangle)]=k-p+1 &\cr
& \hskip 12mm {\rm for} \quad
p\in \{1,3,\ldots,k-1\} & (3.43) \cr
}
$$
$$
\eqalignno{
& dim[W(\vert p;0\rangle,\vert p;3k-p\rangle)]=
dim[W(\vert p+2k;0\rangle,\vert p+2k;3k-p\rangle)]=3k-p+1 &\cr
& \hskip 12mm {\rm for} \quad
p\in \{k+3,k+5,\ldots,2k-1\}. & (3.44) \cr
}
$$
Also in this case the modules with the same dimension cannot be
separated by the Casimir operator. They have however different
spectrum of the Cartan generator $K$. Hence they carry
inequivalent representations of $pB_q$.

\vskip 16pt
\noindent
{\bf 3.2 The case k=odd, m=odd}
\bigskip
\noindent
For any $N\in \Z_+$ (see (2.9))
$a^- \vert p;kN\rangle=0 .$
Therefore the Fock spaces contain more singular vectors than in
the general case (see (3.3)). Now
$$
\vert p;0\rangle,\vert p;k\rangle, \vert p;2k\rangle,\ldots,
\vert p;kN\rangle, \ldots   \eqno(3.45)
$$
are singular vectors in $F(p)$. For each $N$ the subspace
$$
V_{\vert p;kN\rangle}=span\{\vert p;n\rangle \vert n\geq kN\}
\eqno(3.46)
$$
\noindent
is an infinite-dimensional indecomposible  invariant subspace
of $F(p)$. Because of the inclusions
$$
F(p)=V_{\vert p;0\rangle}\supset V_{\vert p;k\rangle}\supset
V_{\vert p;2k\rangle}\supset
\ldots \supset V_{\vert p;kN\rangle} \supset \ldots
$$
one can build up various indecomposible finite-dimensional $pB_q$
modules. For each \hfill \break
$N<M \in {\bf Z}_+$
$$
W(\vert p;kN\rangle,\vert p;kM-1\rangle)
=V_{\vert p;kN\rangle}/V_{\vert p;kM\rangle}
\eqno(3.47)
$$
\noindent
is  an  $(M-N)k$-dimensional indecomposible $pB_q$ module with
singular vectors $\vert p;kN \>, \hfill\break
\vert p; (k+1)N\>,\ldots,\vert
p;k(M-1)\> $.  As before we do not distinguish between the
equivalence classes and their representatives. Since
$$
a^+\xi_{\vert p;kM-1\>}=\xi_{a^+\vert p;kM-1\>}=0 \quad
\Leftrightarrow \quad a^+\vert p;kM-1\>=0 \eqno(3.48)
$$
the transformation of the factor space $(3.47)$ is given with eqs.
(3.11), where one has to replace $2kM-1$ by $kM-1$ in (3.11{\it b,c}).
If $M-N=1 \quad W(\vert p;kN\rangle,\vert p;k(N+1)-1\rangle)
=V_{\vert p;kN\rangle}/V_{\vert p;k(N+1)\rangle}$ is either
irreducible or indecomposible. We are mainly concerned with the
classification of the irreducible $pB_q$ modules. Therefore, as a
first step, we identify some equivalent modules. From Proposition
2 we know that for a given $pB_q$ algebra (i.e. for a fixed
admissible fraction $m\over k$) all irreps can be extracted from
the collection of the modules
$$
W(\vert p;0\rangle,\vert p;2k-1\rangle) \quad 0<Re(p)\leq 4k.
\eqno(3.49)
$$
Accoding to (3.47) (when $N=0$, $M=2$) the above module is
indecomposible. It contains at least two singular vectors, namely
$\vert p;0\>$ and $\vert p;k\>$. The subspace
$W(\vert p;k \rangle,\vert p;2k-1\rangle)$ is an invariant $pB_q$
subspace and therefore the factor space
$$
W(\vert p;0\rangle,\vert p;k-1\rangle)=
W(\vert p;0\rangle,\vert p;2k-1\rangle)/
W(\vert p;k\rangle,\vert p;2k-1\rangle) \eqno(3.50)
$$
is also a $pB_q$ module.

\noindent
\smallskip
\noindent
{\it Proposition 6.} The collection of all $pB_q$ modules
$W(\vert p;0\rangle,\vert p;k-1\rangle)$ is equivalent to the
collection of all invariant subspaces
$W(\vert p;k\rangle,\vert p;2k-1\rangle)$ when $0<Re(p)\leq 4k $.
More precisely,
$$
W(\vert p;0\rangle,\vert p;k-1\rangle)=
W(\vert p+2k;k\rangle,\vert p+2k;2k-1\rangle) \quad
0<Re(p)\leq 2k
\eqno(3.51a)
$$
$$
W(\vert p;0\rangle,\vert p;k-1\rangle)=
W(\vert p-2k;k\rangle,\vert p-2k;2k-1\rangle) \quad
2k<Re(p)\leq 4k.
\eqno(3.51b)
$$
\smallskip
\noindent
{\it Proof.} The transformations of
$W(\vert p;0\rangle,\vert p;k-1\rangle)$ are given with eqs.
(3.19) for $L=k-1$. The same equations describe the action of the
$pB_q$ generators on
$W(\vert p;k\rangle,\vert p;2k-1\rangle)$
for $L=2k-1$ and $n=k,k+1,\ldots,2k-1$. The corresponding
intertwining operator $\varphi$, defined on the basisses, read:
$$
\varphi(\vert p;n\rangle)=\vert p+2\xi k;n+k \rangle
\quad n=0,1,2,\ldots,k-1 \eqno(3.52)
$$
where $\xi =1$ corresponds to $(3.51a)$ and $\xi =-1$
--- to $(3.51b)$.

\smallskip
In view of Proposition 6 we shall consider only the
vacuum modules $W(\vert p;0\rangle,\vert p;k-1\rangle)$
restricting as before the values of $p$ to the interval $(3.16)$.

If $p$ is not an integer the coefficients $\{n+p-1\}$ in $(2.9c)$
and $[n+p-1]$ in $(2.9d)$ never vanish. Therefore the vectors
$(3.45)$ are the only singular vectors in $F(p)$ and the vacuum
$\vert p;0\rangle$ is the only singular vector in
$W(\vert p;0\rangle,\vert p;k-1\rangle)$.
We have obtained the following result.

\smallskip
\noindent
{\it Proposition 7.} The $pB_q$ modules
$W(\vert p;0\rangle,\vert p;k-1\rangle)$
are simple  for $p\notin \{1,2,\ldots,4k\}$. All of them are
$k-$dimensional.
The irreps corresponding to different
$p$ from (3.16) are inequivalent.
The transformation of the basis
$$
\vert p;0\rangle,\;\vert p;1\rangle,\ldots, \vert p;k-1\rangle
\eqno(3.53)
$$
is given with eqs. $(3.19)$ for {\it L=k-1.}

\vskip 14pt
\noindent
{\it 3.2a Representations with even $p$ }
\bigskip
\noindent

{\it Proposition 8.} To each $p\in \{2,4,\ldots,4k\}$ there
corresponds an irreducible $pB_q $  module \hfill\break
$W(\vert p;0\rangle,\vert p;L\rangle)$   with a basis $\vert
p;0\rangle,\vert p;1\rangle,\ldots,\vert p;L\rangle$ and values
of $L$ as follows:
$$
\eqalignno
{
& L=k-p \quad \;\; {\rm for}
\quad  p\in\{2,4,\ldots,k-1\} & (3.54a)\cr
& L=2k-p \quad {\rm for}
\quad p\in\{k+3,k+5,\ldots,2k \} & (3.54b)\cr
& L=3k-p \quad {\rm for}
\quad p\in\{2k+2, 2k+4,\ldots,3k-1\} & (3.54c) \cr
& L=4k-p \quad {\rm for}
\quad p\in\{3k+3, 3k+5,\ldots,4k \} & (3.54d) \cr
& L=k-1 \;\;\quad {\rm for}
\quad p=k+1,\;3k+1. &  (3.54e)\cr
}
$$
The transformations of the basis are given with the  eqs.
(3.19) for the corresponding values of $L$ in (3.54). All
$2k$ modules (3.54) carry inequivalent irreps of $pB_q$.
\smallskip
\smallskip
\noindent
{\it Proof.} We sketch the proof. The modules
$$
W(\vert p;0\rangle,\vert p;k-1\rangle)  \quad {\rm with}
\quad p=k+1,\;3k+1   \eqno(3.55)
$$
remain irreducible, whereas
$$
W(\vert p;0\rangle,\vert p;k-1\rangle) \quad p\in\{2,4,\ldots,4k \}
\;\;p\neq k+1,\;3k+1 \eqno(3.56)
$$
are indecomposible. Each module in (3.56) contains apart from the
vacuum only one more singular vector. Therefore the invariant
subspaces and the factor spaces are irreducible. Each
invariant subspace is equivalent to a factor space:
$$
\eqalignno{
& W(\vert p;k-p+1\rangle,\vert p;k-1\rangle)=
W(\vert p';0\rangle,\vert p';2k-p'\rangle)&\cr
& \hskip 12mm {\rm for} \quad  p=2k-p'+2\in
\{2,4,\ldots,k-1\} & (3.57) \cr
}
$$
$$
\eqalignno{
& W(\vert p;2k-p+1\rangle,\vert p;k-1\rangle)=
W(\vert p';0\rangle,\vert p';3k-p'\rangle)&\cr
& \hskip 12mm {\rm for} \quad  p=4k-p'+2\in
\{k+3,k+5,\ldots,2k \} & (3.58) \cr
}
$$
$$
\eqalignno{
& W(\vert p;3k-p+1\rangle,\vert p;k-1\rangle)=
W(\vert p';0\rangle,\vert p';4k-p'\rangle)&\cr
& \hskip 12mm {\rm for} \quad  p=6k-p'+2\in
\{2k+2,2k+4,\ldots,3k-1\} & (3.59) \cr
}
$$
$$
\eqalignno{
& W(\vert p;4k-p+1\rangle,\vert p;k-1\rangle)=
W(\vert p';0\rangle,\vert p';k-p'\rangle)&\cr
& \hskip 12mm {\rm for} \quad  p=4k-p'+2\in
\{3k+3,3k+5,\ldots,4k \}. & (3.60) \cr
}
$$
Thus, (up to equivalence) one can consider only the vacuum
modules (3.54).  Some of them have one and the same dimension,
namely
$$
\eqalignno{
& dim[W(\vert p;0\rangle,\vert p;k-p\rangle)]=
dim[W(\vert p+2k;0\rangle,\vert p+2k;k-p\rangle)]=k-p+1 &\cr
& \hskip 12mm {\rm for} \quad
p\in \{2,4,\ldots,k-1\} & (3.61) \cr
}
$$
$$
\eqalignno{
& dim[W(\vert p;0\rangle,\vert p;2k-p\rangle)]=
dim[W(\vert p+2k;0\rangle,\vert p+2k;2k-p\rangle)]=2k-p+1 &\cr
& \hskip 12mm {\rm for} \quad
p\in \{k+3,k+5,\ldots,2k\} & (3.62) \cr
}
$$
$$
dim[W(\vert k+1;0\rangle,\vert k+1;k-1\rangle)]=
dim[W(\vert 3k+1;0\rangle,\vert 3k+1;k-1\rangle)]=k.\eqno(3.63)
$$

The modules with the same dimension cannot be separated by the
Casimir operator.  The spectrum of $K$ is however different.
Therefore these modules are inequivalent.

\vskip 14pt
\noindent
{\it 3.2b Representations with odd $p$ }
\bigskip
\noindent
This collection of representations is somewhat different from the
other cases.  All Fock modules
$W(\vert p;0\rangle,\vert p;k-1\rangle)$
remain irreducible for $p\in \{1,3,5,\ldots,4k-1\}$.

\vfill \eject
\noindent
{\bf 3.3 The case k=odd, m=even}
\bigskip
\noindent

According to {\it Proposition 2} the irreps (up to equivalence)
are realized in the vacuum modules $W(\vert p;0\rangle,\vert
p;2k-1\rangle)$ with $0<Re(p)\leq 4k$. For the algebras from this
class one can further restrict the values of $p$.

\smallskip
\noindent
{\it Proposition 9}. The following modules are equivalent:
$$
W(\vert p;0\rangle,\vert p;2k-1\rangle)=W(\vert
p+2k;0\rangle,\vert p+2k;2k-1\rangle) \quad {\rm for} \quad 0<Re(p)\leq
2k \eqno(3.64)
$$
$$
\hskip 10pt W(\vert p;0\rangle,\vert p;2k-1\rangle)=W(\vert
p+k;0\rangle,\vert p+k;2k-1\rangle) \;\; {\rm for} \;\; 0<Re(p)\leq k
\;\;
{\rm and} \;\; m=4(mod\;4). \eqno(3.65)
$$
\smallskip
\noindent
{\it Proof.} The intertwining operators, written on the basisses read
$$
\varphi (\vert p;n\rangle)=\vert p+\alpha k;n\rangle \quad n=0,1,
\ldots ,2k-1  \eqno(3.66)
$$
where $\alpha =2$ for (3.64)  and $\alpha =1$ for (3.65).

\smallskip
Hence, without loss of generality we assume
$$
0<Re(p)\leq 2k \quad {\rm if} \quad m=2(mod\;4) \eqno(3.67)
$$
$$
0<Re(p)\leq k \quad {\rm if} \quad m=4(mod\;4). \eqno(3.68)
$$

\noindent
{\it Proposition 10.} The $pB_q$ module $W(\vert p;0\rangle,\vert
p;2k-1\rangle)$ is irreducible if $p$ is not an integer. All such
modules are $2k$-dimensional.  The irreps corresponding to
different $p$ from (3.67) and (3.68) are inequivalent.  The
transformations of the basis (3.15) is described with eqs. (3.19)
for $L=2k-1$.

\smallskip
\noindent
{\it Proof.} The same as in  Proposition 3.

If $p$ is an integer each space $W(\vert p;0\rangle,\vert
p;2k-1\rangle)$ is indecomposible. Apart from $\vert p;0\rangle$
it contains only one more singular vector $\vert p;L+1\rangle$, where
$$
\eqalignno{
&L=2k-p \quad {\rm for} \quad p={\rm even} & (3.69a) \cr
&L=k-p \quad \;\; {\rm for} \quad p\in \{1,3,\ldots ,k\} & (3.69b) \cr
&L=3k-p \quad {\rm for} \quad p\in \{k+2,k+4,\ldots ,2k-1\} \quad
(m=2(mod\;4)). & (3.69c) \cr
}
$$
In all cases $W(\vert p;L+1\rangle,\vert p;2k-1\rangle)$ is an
irreducible $pB_q$ module. So is the corresponding factor space
$$
W(\vert p;0\rangle,\vert p;L\rangle)=W(\vert p;0\rangle,\vert
p;2k-1\rangle)/W(\vert p;L+1\rangle,\vert p;2k-1\rangle). \eqno(3.70)
$$
\noindent
Each invariant subspace is equivalent to a certain factor space:
$$
\eqalignno
{
& W(\vert 1;k\rangle,\vert 1;2k-1\rangle)=
W(\vert 1;0\rangle,\vert 1;k-1\rangle) &  (3.71) \cr
& W(\vert p;k-p+1\rangle,\vert p;2k-1\rangle)=
  W(\vert p';0\rangle,\vert p';3k-p'\rangle) & \cr
& \hskip 12mm {\rm for} \quad m=2(mod\;4), \;\;
  p=2k-p'+2\in \{3,5,\ldots ,k\} & (3.72) \cr
& W(\vert p;3k-p+1\rangle,\vert p;2k-1\rangle)=W(\vert
  p';0\rangle,\vert p';k-p'\rangle) & \cr
& \hskip 12mm {\rm for} \quad m=2(mod\;4), \;\; p=2k-p'+2\in
  \{k+2,k+4,\ldots ,2k-1\}& (3.73) \cr
& W(\vert p;k-p+1\rangle,\vert p;2k-1\rangle)=W(\vert
  p';0\rangle,\vert p';2k-p'\rangle)& \cr
& \hskip 12mm {\rm for} \quad m=4(mod\;4), \;\; p=k-p'+2\in
  \{3,5,\ldots ,k\}& (3.74) \cr
& W(\vert p;2k-p+1\rangle,\vert p;2k-1\rangle)=
  W(\vert p';0\rangle,\vert p';2k-p'\rangle)& \cr
& \hskip 12mm {\rm for} \quad  p=2k-p'+2={\rm even}
& (3.75) \cr
}
$$
\noindent
We skip the explicit form of the intertwining operators and
collect the result in  a proposition.

\smallskip
\noindent
{\it Proposition 11.} To each integer $p$, $0<p\leq 2k$ if
$m=2(mod \; 4)$, and $0<p\leq k$ if $m=4(mod \; 4)$,
there corresponds an irreducible $pB_q$ module
$W(\vert p;0\rangle,\vert p;L\rangle) $ with a basis
$\vert p;0\rangle,\vert p;1\rangle,\ldots ,\vert p;L\rangle$. All
such modules are inequivalent. Their transformations under the
action of $pB_q$ are given with eqs. (3.19) for the corresponding
values of $L$ in (3.69).

\vskip 32pt
\noindent
{\bf 4. Unitarizable representations}

\bigskip
\noindent
In the present Section we classify the unitarizable Fock
representations of the deformed para-Bose superalgebra $pB_q$.
The concept of an unitarizable representation of an arbitrary
associative algebra $A$ depends on the definition of the
antilinear antiinvolution $\omega : A \rightarrow A$ and on the
metric in the corresponding $A$-module. One and the same
representation can be unitarizable with respect to one
antiinvolution and not unitarizable with respect to another one.

Having in mind the physical condition (1.1) and the requirement
the "Hamiltonian" $H$ to be Hermitian operator, we define
$\omega$ on the generators $a^\pm$ and $K$ as
$$
\omega (a^\pm)=a^\mp \quad\quad \omega (K^{\pm 1})=K^{\mp 1}
\eqno(4.1)
$$
and extend it on $pB_q$ as an antilinear antiinvolution,
namely
$$
\omega(\alpha a + \beta b)=\alpha^* \omega(a)+\beta^*\omega(b)
\quad \omega(ab)=\omega(b)\omega(a) \eqno(4.2)
$$
for all $a,b \in pB_q$ and $\alpha, \beta \in \C$
with $\alpha^*$ being the complex conjugate of $\alpha$.

The representation of $pB_q$ in a Hilbert space $W$ with scalar
product $(\;,\;)$ is unitarizable if
$$
(ax,y)=(x,\omega(a)y) \quad \forall\;a\in pB_q \;\;{\rm and}
\;\; x,y\in W \eqno(4.3)
$$
or, equivalently, if
$\omega(a)=a^\dagger $. On the generators of $pB_q$ the latter
condition reads:
$$
(a^-)^\dagger=a^+ \quad K^\dagger=K^{-1}. \eqno(4.4)
$$

The problem is to select those irreducible modules
$W(\vert p;0 \rangle ,\vert p;L \rangle)$ for which the unitarity
condition (4.4) can be satisfied. To this end we introduce a new basis
$$
\vert p;n \rangle=\alpha(p;n)\vert p;n) \quad n=1,2,\ldots,L
\quad \alpha(p;n)\in \C \eqno(4.5)
$$
which is declared to be orthonormed. Then the unitarity condition
is equivalent to the requirement the following two conditions to
be satisfied:
$$
 \left|{{\alpha(p;n)}\over{\alpha(p;n+1)}}  \right|^2
  ={{{2sin({\pi\over 2}{m\over k}(n+p))}cos({\pi\over 2}{m\over k}(n+1))}
  \over {sin({{\pi m}\over k}})} \quad\quad {\rm n=even} \eqno(4.6a)
$$
$$
 \left|{{\alpha(p;n)}\over{\alpha(p;n+1)}}  \right|^2
  ={{{2sin({\pi\over 2}{m\over k}(n+1))}cos({\pi\over 2}{m\over k}(n+p))}
  \over {sin({{\pi m}\over k}})} \quad\quad {\rm n=odd}. \eqno(4.6b)
$$

\noindent
The unknowns in the above equations are the algebras $pB_q$, i.e.,
the admissible pairs $m\over k$ and the irreducible modules
$W(\vert p;0 \rangle ,\vert p;L \rangle)$, i.e., the values of
$p$ and $L$. The eqs. (4.6) have solutions only if the r.h.s. of
them is nonnegative number for any $n$ from the basis in
$W(\vert p;0 \rangle ,\vert p;L \rangle)$. Thus, the problem
is to solve a set of inequalities. Bellow we list the algebras
$pB_q$ with $q=e^{i{\pi \over 2}{m\over k}}$ in terms of the
admissible $m$ and $k$ and their unitarizable representations.

$$
\vcenter{\openup3\jot\halign{#\hfil & \hskip 12pt # \hfil \cr
The algebra $pB_q$ & Unitarizable modules \cr
(1)$\;\;$   $ m=1,\; k=3,5,7,\ldots$
& $W(\vert p;0 \rangle ,\vert p;k-1 \rangle) \quad
0 \leq p \leq 2$  \cr
(2)$\;\;$   $m=1,\;k=2,4,6,\ldots $
& $W(\vert p;0 \rangle ,\vert p;k-p \rangle) \quad
 p=1,3,5,\ldots,k-1 $  \cr
(3)$\;\;$   $m=1,\; k=3,5,7,\ldots $
& $W(\vert p;0 \rangle ,\vert p;k-p \rangle) \quad
 p=2,4,6 \ldots,k-1 $  \cr
 (4)$\;\;$   $m=1(mod\; 4),\; k=2,3,4,\ldots $
& $W(\vert k-1;0 \rangle ,\vert k-1;1 \rangle)$  \cr
 (5)$\;\;$   $m=3(mod\; 4),\; k=2,3,4,\ldots $
& $W(\vert 3k-1;0 \rangle ,\vert 3k-1;1 \rangle)$  \cr
 (6)$\;\;$   $m=3, \; k=10,12,14,\ldots $
& $W(\vert 3k-3;0 \rangle ,\vert 3k-3;3 \rangle).$  \cr
}} \eqno(4.7)
$$
\smallskip
\noindent
We underline that everywhere in eqs. (4.7) only the admissible
values of $m$ and $k$ are to be considered (see (3.2)).

The above equations indicate that the algebras $pB_q$
corresponding to $ m=1$ and any odd $k$ have a continuous class
of unitarizable representations. In all other cases
the number of the unitarizable irreps is finite and in fact each
algebra with $m\neq 1$ has no more than two representation. In
the cases (4) and (5) the representation is 2-dimensional. In the
new basis (4.5) the matrices of the generators read:

$$ a^-=\left(\matrix{ 0 &
\sqrt{{cos{\pi\over 2}{m\over k}}\over{sin{\pi\over 2}{m\over
k}}}\cr 0 & 0 \cr }\right) \quad a^+=\left(\matrix{ 0 & 0 \cr
\sqrt{{cos{\pi\over 2}{m\over k}}\over{sin{\pi\over 2}{m\over
k}}} & 0 \cr
}\right) \quad
K=\left(\matrix{
ie^{-i{\pi\over 2}{m\over k}} & 0 \cr
0 & ie^{i{\pi\over 2}{m\over k}}  \cr
}\right) \eqno(4.8)
$$

\noindent
Similarly, the 4-dimensional representation from the case (6)
reads:
$$
a^-=\left(\matrix{
0 & \sqrt{{cos{9\pi\over {2k}}}\over {sin{3\pi\over {2k}}}} & 0 &
0 \cr
0 & 0 & \sqrt{2sin{3\pi\over k}} & 0 \cr
0 & 0 & 0 & \sqrt{{cos{9\pi\over {2k}}}\over {sin{3\pi\over
{2k}}}} \cr
0 & 0 & 0 & 0 \cr
}\right) \quad
K=\left(\matrix{
ie^{-{9i\pi\over {2k}}} & 0 & 0 & 0 \cr
0 & ie^{-{3i\pi\over {2k}}} & 0 & 0 \cr
0 & 0 & ie^{3i\pi\over {2k}} & 0  \cr
0 & 0 & 0 & ie^{9i\pi\over {2k}} \cr
}\right) \eqno(4.9)
$$
with $a^+$ represented by the transposed matrix of $a^-$. Note
that the algebras with $m=3$ and $k=10,12,14,\ldots$ have only
two unitarizable irreps, namely (4.8) and (4.9); the algebras
with $m=3$ and $k=2,4,6,8$
and those with $m=5,7,9,\ldots$  have only a 2-dimensional
unitarizable irrep.  The algebras $pB_q$ with even $m$ have no
unitarizable representations at all.

The transformation relations of all unitarizable modules can be
written in a compact form. In the orthonormed basis
$\vert p;n)$ eqs. (3.19) read:
$$
\eqalignno
{
& K\vert p;n)=e^{i{\pi \over 2}{m\over k}(2n+p)}\vert p;n) &
  (4.10a) \cr
&&\cr
& a^-\vert p;n)=\sqrt
  {{{2sin({\pi\over 2}{m\over k}n)}cos({\pi\over 2}{m\over k}(p+n-1))}
  \over {sin({{\pi m}\over k}})}\vert p;n-1) \quad\quad {\rm n=even} &
  (4.10b) \cr
&&\cr
& a^-\vert p;n)=\sqrt
  {{{2sin({\pi\over 2}{m\over k}(p+n-1))}cos({\pi\over 2}{m\over k}n)}
  \over {sin({{\pi m}\over k}})}\vert p;n-1) \quad\quad {\rm n=odd} &
  (4.10c) \cr
&&\cr
& a^+\vert p;n)=\sqrt
  {{{2sin({\pi\over 2}{m\over k}(p+n))}cos({\pi\over 2}{m\over k}(n+1))}
  \over {sin({{\pi m}\over k}})}\vert p;n+1) \quad\quad {\rm n=even} &
  (4.10d) \cr
&&\cr
& a^+\vert p;n)=\sqrt
  {{{2sin({\pi\over 2}{m\over k}(n+1))}cos({\pi\over 2}{m\over k}(p+n))}
  \over {sin({{\pi m}\over k}})}\vert p;n+1) \quad\quad {\rm n=odd} &
  (4.10e) \cr
}
$$
We have skipped eq. (3.19{\it f}), since it is automatially
satisfied: the creation operator $a^+$, which is the negative
root vector, annihilates the lowest weight vectors $\vert p;L)$
within each unitarizable module.

\vskip 32pt
\noindent
{\bf 5. Concluding remarks and discussions}

\bigskip
\noindent
We have studied  root of unity representations of the deformed
para-Bose algebra $pB_q=U_q[osp(1/2)]$ with a particular emphasis
on the unitarizable irreps. All of them are realized in
finite-dimensional modules with a highest and a lowest weight.
The irreps from Sec. 3.2 and also all irreps, corresponding to
integer $p$ (except $p=k+1,\;3k+1$ in Sec. 3.1{\it b}) are new.

In the nondeformed case the representations of the para-Bose
operators, corresponding to an order of the statistics $p=1$
reduce to usual Bose operators [42]. In [1] it was shown
that similar relation holds  in the deformed case for generic
$q$. It is straightforward to check that in the cases $p=1$,
$m=1$, $k=2,3,\ldots$ eqs. (4.10) recover also all root of unity
unitarizable irreps of the deformed Bose operators [24-27] as
given in [1].

Using the approach of the present paper one can try to construct
representations (including root of 1 representations) for
$pB_q(n)=U_q[osp(1/2n)]$. To this end one can use $n$-pairs of
deformed pB operators as given in [2,4,5]. The solution, however,
is not going to be easy for arbitrary values of $p$, if one takes
into account that the problem has not been solved even in the
nondeformed case. Only the case with $p=1$ is easy. It leads
directly to root of 1 representations of $U_q[osp(1/2n)]$, if one
uses $q$-commuting deformed Bose operators as defined in [5].
Other root of 1 representations based on a realization with
commuting $q$-Bose operators (which means also the case $p=1$)
were obtained in [44]. In this relation we note that $n$ pairs of
commuting deformed Bose operators are already generators of
$U_q[osp(1/2n)]$ (in the $q$-Bose representation). Therefore they
provide the simplest $q$-Boson realization of $U_q[osp(1/2n)]$
[28].

Finally we mention that all our representations correspond to
$q$ being even root of unity: $q^{4k}=1$.  In case of deformed
simple Lie algebras this seems to be the more difficult case.
Complete results exist for $q$ being only odd roots of 1
[45].
\bigskip
\noindent
{\bf Acknowledgments}

\smallskip
\noindent
The authors would like to thank Prof. Randjbar-Daemi for the kind
hospitality at the High Energy Section of ICTP, Trieste.
Constructive discussions with Prof. J. Van der Jeugt and Prof. R.
Jagannathan are greatly acknowledged. The work was supported by
the Grant $\Phi$-416 of the Bulgarian Foundation for Scientific
Research.


\vfill\eject

{\bf References}

\vskip 12pt
{\settabs \+  [11] & I. Patera, T. D. Palev, Theoretical
   interpretation of the experiments on the elastic \cr

\+ [1] & Celeghini E, Palev T D and Tarlini M 1990 {\it Preprint}
         YITP/K-865 Kyoto and \cr
\+     & 1991  {\it Mod. Phys. Lett. B} {\bf 5} 187 	\cr

\+ [2] & Palev T D 1993 {\it J. Phys. A} {\bf 26} L1111  \cr

\+ [3] & Hadjiivanov L K 1993 {\it Journ. Math. Phys.} {\bf 34}
         5476 \cr

\+ [4] & Palev T D 1994  Quantization of the Lie algebra $so(2n+1)$
         and of the Lie superalgebra  \cr
\+     & $osp(1/2n)$ with preoscillator generators {\it Preprint}
         TWI-94-29  University of Ghent \cr
\+     & and {\it Journ. Group Theory in Physics} {\bf 3})
		 (to appear)\cr

\+ [5] & Palev T D and Van der Jeugt J 1995
         {\it J. Phys. A: Math. Gen.} {\bf 28} 2605  \cr

\+ [6] & Palev T D 1994 {\it Lett. Math. Phys.} {\bf 31} 151 \cr

\+ [7] & Greenberg O W and Mohapatra R N 1987 {\it Phys.\
         Rev.\ Lett.} {\bf 59} 2507 \cr

\+ [8] & Floreanini R and Vinet L 1990 {\it J.\ Phys.\ A~: Math.\
         Gen.} {\bf 23} L1019 \cr
\+ [9] & Odaka K, Kishi T and Kamefuchi S 1991
	 {\it J.\ Phys.\ A~: Math.\ Gen.} {\bf 24} L591 \cr

\+ [10] & Beckers J and Debergh N 1991
          {\it J.\ Phys.\ A~: Math.\ Gen.} {\bf 247} L1277 \cr

\+ [11] & Chaturvedi S and Srinivasan V 1991 {\it Phys.\ Rev.\ A}
          {\bf 44} 8024  \cr

\+ [12] & Krishna-Kumari M Shanta P Chaturvedi S and Srinivasan V\cr
\+      & 1992 {\it Mod.\ Phys.\ Lett.\ A } {\bf 7 } 2593 \cr

\+ [13] & Bonatsos D and Daskaloyannis C 1993 {\it
          Phys.\ Lett.\ B} {\bf 307} 100 and the references therein \cr

\+ [14] & Flato M, Hadjiivanov L K and Todorov I T 1993
         {\it Found.\ Phys.} {\bf 23 } 571 \cr

\+ [15] & Macfarlane A J 1993 Generalized oscillator systems and
          and their  parabosonic  \cr
\+      & interpretation {\it Preprint}  DAMPT 93-37 \cr

\+ [16] & Palev T D 1993 {\it J.\ Math.\ Phys.} {\bf 34} 4872 \cr

\+ [17] & Quesne C 1994 {\it Phys.\ Lett.\ A} {\bf 193} 245\cr

\+ [18] & Van der Jeugt J and Jagannathan R 1994  Polynomial
          deformations of $osp(1/2)$ and  \cr
\+  	& generalized parabosons hep-th/9410145 and
          {\it J.\ Phys.\ A~: Math.\ Gen.} (to appear) \cr

\+ [19] & Macfarlane A J 1994 {\it  J.\ Math.\ Phys.} {\bf 35}
          1054 \cr

\+ [20] & Cho K H, Chaiho Rim, Soh D S and Park S U 1994
          {\it J.\ Phys.\ A~: Math.\ Gen.} {\bf 27} 2811 \cr

\+ [21] & Chakrabarti R and Jagannathan R 1994
         {\it J.\ Phys.\ A~: Math.\ Gen.} {\bf 27} L277; \cr
\+      & 1994 {\it Int. J. Mod. Phys. A } {\bf 9} 1411 \cr

\+ [22] & Green H S 1994 {\it Austr.\ J.\ Phys.} {\bf 47} 109 \cr

\+ [23] & Ganchev A and Palev T D 1978 {\it Preprint} JINR P2-11941;
          1980 {\it J. Math. Phys.} {\bf 21} 797 \cr

\+ [24] & Macfarlane A J 1989 {\it J.\ Phys.\ A~: Math.\ Gen.}
          {\bf 22}  4581  \cr

\+ [25] & Biedenharn L C 1989 {\it J.\ Phys.\ A~: Math.\ Gen.} {\bf 22}
          L873 \cr

\+ [26] & Sun C P and Fu H C 1989 {\it J.\ Phys.\ A~: Math.\ Gen.}
          {\bf 22}  L983  \cr

\+ [27] & Hayashi T 1990 {\it Commun.\ Math.\ Phys.} {\bf 127}  129  \cr

\+ [28] & Palev T D 1993 {\it Lett. Math. Phys.} {\bf 28} 187 \cr

\+ [29] & Palev T D 1994 {\it J. Phys. A} {\bf 27} 7373  \cr

\+ [30] & Chaichan M and Kulish P 1990 {\it Phys.\ Lett.\ B}
          {\bf 234} 72 \cr

\+ [31] & Bracken A J, Gould M D and Zhang R B 1990
          {\it Mod. Phys. Lett. A} {\bf 5} 331 \cr

\+ [32] & Floreanini R, Spiridonov V P and Vinet L 1990
          {\it Phys.\ Lett.\ B} {\bf 242} 383 \cr

\+ [33] & Floreanini R, Spiridonov V P and Vinet L 1991
          {\it Commun.\ Math.\ Phys.}{\bf 137} 149 \cr

\+ [34] & d\'{}$\!\!$ Hoker, Floreanini R and Vinet L 1991
          {\it Journ. Math. Phys.} {\bf 32} 1247 \cr

\+ [35] & Khoroshkin S M  and Tolstoy  V N  1991 {\it Commun.\
          Math.\ Phys.} {\bf 141} 599 \cr

\+ [36] & Kulish P P and Reshetikhin N Yu 1989 {\it Lett. Math. Phys}
          {\bf 18} 143 \cr

\+ [37] & Saleur H 1990 {\it Nucl. Phys. B} {\bf 336} 363 \cr

\+ [38] & Sun Chang-Pu, Fu Hong-Chen and Ge Mo-Lin 1991 {\it
          Lett. Math. Phys} {\bf 23} 19 \cr

\+ [39] & Ge Mo-Lin, Sun Chang-Pu and Xue Kang 1992 {\it Phys.
          Lett. A} {bf 163} 176 \cr

\+ [40] & Sun Chang-Pu 1993 {\it N. Cim. B} {\bf 108} 499 \cr

\+ [41] & Ubriaco M R 1993 {\it Mod. Phys. Lett. A} {\bf 8} 89 \cr

\+ [42] & Green H S 1953 {\it Phys. Rev.} {\bf 90} 270\cr

\+ [43] & Pasquier V and Saleur H 1990 {\it Nucl. Phys. B}
          {\bf 330} 523 \cr

\+ [44] & Fu Hong-Chen and Ge Mo-Lin 1992 {\it Commun. Theor. Phys.}
          {\bf 18} 373 \cr

\+ [45] & De Concini C and Kac V G 1990 Colloque Dixmier
          {\it Progr. Math.} {\bf 92} 471  \cr

\end